\documentclass[10pt,twoside,reqno]{amsart}
\usepackage{amssymb}
\usepackage{multirow}
\calclayout

\numberwithin{equation}{section}
\makeatletter

\renewcommand{\@secnumfont}{\bfseries}

\renewcommand{\section}{\@startsection{section}{1}%
  {0mm}{.7\linespacing\@plus\linespacing}{.5\linespacing}
  {\normalfont\bfseries\centering}}

\newcommand{\bibsection}{\@startsection{section}{1}%
  {0mm}{.7\linespacing\@plus\linespacing}{.5\linespacing}
  {\normalfont\scshape\centering}}

\renewcommand{\@biblabel}[1]{#1.}

\newtheorem{thm}{\bf Theorem}[section]

\begin{document}

\begin{center}
	
	{\large\bf Degenerate central factorial numbers of the second kind
		\rule{0mm}{6mm}\renewcommand{\thefootnote}{}

		\footnotetext{
			\rule{2.4mm}{0mm}Keywords and Phrases. degenerate central factorial numbers of the second kind, degenerate central factorial polynomials of the second kind. \\
			
	}}
	
	\vspace{1cc}
	{\large\it Taekyun Kim, Dae San Kim}
	
	\vspace{1cc}
	\parbox{24cc}{{\small
			In this paper, we introduce the degenerate central factorial polynomials and numbers of the second kind which are degenerate versions of the central factorial polynomials and numbers of the second kind. We derive some properties and identities for those polynomials and numbers. We obtain, among other things, recursive formulas for the degenerate central factorial polynomials and numbers of the second kind.
			
			%
			
	}}
\end{center}

\vspace{1cc}


\vspace{1.5cc}
\begin{center}
	{\bf 1. Introduction}
\end{center}

Various degenerate versions of special polynomials and numbers have drawn the attention of many mathematicians in recent years. The origin of these are the papers by Carlitz \cite{bib2,bib3} on degenerate Bernoulli and degenerate Euler polynomials and numbers. The degenerate Bernoulli polynomials were later rediscovered by Ustinov \cite{bib13} under the name of Korobov polynomoals of the second. Also, Korobov \cite{bib11} introduced Korobov polynomials of the first kind which are in fact a degenerate version of Bernoulli polynomials of the second kind. All of them studied some arithmetic and combinatorial aspects of those degenerate special polynomials and numbers.\\
\indent More recently, along the same line, studying various degenerate versions of many special polynomials and numbers regained attention of the present authors, their colleagues and some other people in connection with their interest not only in arithmetic and combinatorial properties but also in certain symmetric identities and differential equations \cite{bib9,bib10}. This idea of introducing some degenerate version of certain polynomials and numbers has been extended even to transcendental functions so that degenerate gamma functions were introduced in \cite{bib7}.\\
\indent Here in this paper we study the degenerate central factorial polynomials and numbers of the second kind which are degenerate versions of the central factorial polynomials and numbers of the second kind. We derive some properties and identities for those polynomials and numbers. In particular, we will be able to find recursive formulas for the degenerate central factorial polynomials and numbers of the second kind. As to  degenerate central factorial numbers of the first kind, we will be content with defining them. \\
\indent In conclusion, we may say that studying some degenerate versions of certain special polynomials and numbers are promising area of research and that there are still many things yet to be uncovered.

For $ n \in \mathbb{N} \cup \{0\}$, it is known that the Stirling numbers of the first kind are defined as
\begin{equation} \label{01}
(x)_n = x(x-1)(x-2) \cdots (x-n+1) = \sum_{l=0}^n S_1(n,l) x^l, \,(n\geq 1),\,(x)_0 = 1.
\end{equation}

\noindent  As shown in \cite{bib12}, the Stirling numbers of the first kind satify the relation
\begin{equation} \label{02}
S_1(n+1,k) = S_1(n,k-1) - n S_1(n,k), \,\,(1 \leq k \leq n).
\end{equation}

Carlitz \cite{bib2,bib3} studied the degenerate Euler polynomials given by
\begin{equation} \label{03}
\frac{2}{(1+\lambda t)^{\frac{1}{\lambda}} +1} (1+\lambda t)^{\frac{x}{\lambda}} = \sum_{n=0}^\infty \mathcal{E}_{n,\lambda} (x) \frac{t^n}{n!} ,\,\, (\lambda \in \mathbb{R}).
\end{equation}

Now, as a generalization of the falling factorial sequence in \eqref{01}, the $\lambda$-analogue of the falling factorial sequence are defined as follows:
\begin{equation} \label{04}
(x)_{0,\lambda} = 1, (x)_{n,\lambda} = x(x-\lambda)(x-2\lambda) \cdots (x-(n-1)\lambda), (n \geq 1).
\end{equation}

As defined in \cite{bib7}, the $\lambda$-Stirling numbers of the second kind are given by
\begin{equation} \label{05}
\frac{1}{k!} ((1+\lambda t)^{\frac{1}{\lambda}}-1)^{k} = \sum_{n=k}^\infty S_{2,\lambda}(n,k) \frac{t^n}{n!},\, (k \in \mathbb{N}).
\end{equation}

\noindent Note that, by taking the limit $\lambda \rightarrow 0$ in \eqref{05}, we have  $ \lim_{\lambda \rightarrow 0} S_{2,\lambda}(n,k) = S_{2}(n,k), (n, k \geq 0)$.
Here, as we can see in \cite{bib5,bib8,bib12}, $S_{2}(n,k)$ are the Stirling numbers of the second kind given by
\begin{equation*}
	x^n = \sum_{l=0}^n S_2(n,l) (x)_l, (n \geq 0).
\end{equation*}

As shown in \cite{bib8}, the $\lambda$-analogue of binomial expansion is given by
\begin{equation} \label{06}
(1+\lambda t)^{\frac{x}{\lambda}} = \sum_{l=0}^\infty \binom{x}{l}_{\lambda} t^l,
\end{equation}
where $\binom{x}{l}_{\lambda} = \frac{(x)_{l,\lambda}}{l!} = \frac{x(x-\lambda)(x-2\lambda) \cdots (x-(n-1)\lambda)}{l!}, \,(l \geq 1), \binom{x}{0}_{\lambda}=1$

The central factorial $x^{[n]}$ is defined by the generating function
\begin{equation} \label{07}
\sum_{n=0}^\infty x^{[n]} \frac{t^n}{n!} = \bigg(\frac{t}{2} + \sqrt{1+\frac{t^2}{4}} \bigg)^{2x}.
\end{equation}

\noindent From \eqref{07}, we note that
\begin{equation*}
	x^{[n]} = x(x+\frac{n}{2}-1) (x+\frac{n}{2}-2) \cdots (x-\frac{n}{2}+1),\,\,(n \geq 1),\,\,x^{[0]} =1.
\end{equation*}

As defined in \cite{bib1,bib4,bib6,bib14}, for any nonnegative integer $n$,  the central factorial numbers of the first kind are given by
\begin{equation} \label{08}
x^{[n]} = \sum_{k=0}^n  t(n,k) x^k.
\end{equation}

\noindent By \eqref{08}, we easily get
\begin{equation} \label{09}
\frac{1}{k!} \bigg(2 \log \bigg(\frac{t}{2} + \sqrt{1+\frac{t^2}{4}}\bigg) \bigg)^k  = \sum_{n=k}^\infty  t(n,k) \frac{t^n}{n!}.
\end{equation}

Let $f(t) = 2 \log(\frac{t}{2} + \sqrt{1+\frac{t^2}{4}})$. Then we have
\begin{equation} \label{10}
f^{-1}(t) = e^{\frac{t}{2}} - e^{-\frac{t}{2}}.
\end{equation}

\noindent In view of \eqref{09} and \eqref{10}, we define the central factorial numbers of the second kind by
\begin{equation}  \label{11}
\frac{1}{k!} (e^{\frac{t}{2}} - e^{-\frac{t}{2}})^{k} = \sum_{n=k}^\infty T(n,k) \frac{t^n}{n!}.
\end{equation}

\noindent Thus, as shown in \cite{bib1,bib4,bib6,bib14}, and from \eqref{11} we easily get
\begin{equation}  \label{12}
x^n = \sum_{k=0}^n T(n,k) x^{[k]}.
\end{equation}

\noindent From \eqref{12}, we note that
\begin{equation*}
	T(n,k) = T(n-2,k-2) + \frac{k^2}{4} T(n-2,k), (n, k \geq 2).
\end{equation*}

\vspace{1.5cc}
\begin{center}
	{\bf 2. A note on Central factorial numbers and polynomials of the second kind}
\end{center}

\vspace{0.1in}

The central difference operator $\delta$ is defined by
\begin{equation}\label{13}
\delta f(x) = f(x+\frac{1}{2}) - f(x-\frac{1}{2}).
\end{equation}

\noindent By proceeding induction with \eqref{13}, we can easily show that
\begin{equation}\label{14}
\delta^k f(x) = \sum_{l=0}^k \binom{k}{l} f(x+l-\frac{k}{2}) (-1)^{k-l}, (k \in \mathbb{N}).
\end{equation}

\noindent From \eqref{14}, we note that

\begin{equation}\label{15}
\begin{split}
\delta^k x^{m+1} & = \sum_{l=0}^k \binom{k}{l} (x+l-\frac{k}{2})^{m+1} (-1)^{k-l} = \sum_{l=0}^k \binom{k}{l} (x+l-\frac{k}{2})^{m} (-1)^{k-l} (x+l-\frac{k}{2}) \\
& = (x-\frac{k}{2}) \delta^k x^{m} + k \sum_{l=1}^k \binom{k-1}{l-1} (x+l-\frac{k}{2})^{m} (-1)^{k-l} \\
& = (x-\frac{k}{2}) \delta^k x^{m} + k \sum_{l=1}^k \bigg\{\binom{k}{l} - \binom{k-1}{l} \bigg\} (x+l-\frac{k}{2})^{m} (-1)^{k-l} \\
& = (x-\frac{k}{2}) \delta^k x^{m} + k \sum_{l=0}^k \bigg\{\binom{k}{l} - \binom{k-1}{l} \bigg\} (x+l-\frac{k}{2})^{m} (-1)^{k-l} \\
& = (x-\frac{k}{2}) \delta^k x^{m} + k \bigg(\delta^k x^{m} + \delta^{k-1} x^{m} \bigg) = (x+\frac{k}{2}) \delta^k x^{m} + k \delta^{k-1} (x-\frac{1}{2})^{m}. \\
\end{split}
\end{equation}

We define the degenerate central factorial polynomials of the second kind by
\begin{equation}\label{16}
\frac{1}{k!} (1+\lambda t)^{\frac{x}{\lambda}} ((1+\lambda t)^{\frac{1}{2 \lambda}} - (1+\lambda t)^{-\frac{1}{2 \lambda}})^k = \sum_{n=k}^\infty T_{2,\lambda}(n,k \mid x) \frac{t^n}{n!}.
\end{equation}

\noindent When $x=0$, $T_{2,\lambda}(n,k)=T_{2,\lambda}(n,k|0)$ are called the degenerate central factorial numbers of the second kind so that
\begin{equation}\label{17}
\frac{1}{k!}((1+\lambda t)^{\frac{1}{2 \lambda}} - (1+\lambda t)^{-\frac{1}{2 \lambda}})^k = \sum_{n=k}^\infty T_{2,\lambda}(n,k) \frac{t^n}{n!}.
\end{equation}

From \eqref{17} and with the notation in \eqref{04}, we have
\begin{equation}\label{18}
\begin{split}
& \frac{1}{k!} (1+\lambda t)^{\frac{x}{\lambda}} ((1+\lambda t)^{\frac{1}{2 \lambda}} - (1+\lambda t)^{-\frac{1}{2 \lambda}})^k \\
& = \bigg(\sum_{m=0}^{\infty} (x)_{m,\lambda} \frac{t^m}{m!} \bigg) \bigg(\sum_{l=k}^{\infty} T_{2,\lambda}(l,k) \frac{t^l}{l!} \bigg) \\
& = \sum_{n=k}^{\infty} \bigg\{\sum_{l=k}^{n} \binom{n}{l} T_{2,\lambda}(l,k) (x)_{n-l,\lambda} \bigg\} \frac{t^n}{n!}.
\end{split}
\end{equation}

\noindent Therefore, by \eqref{16} and \eqref{18}, we obtain the following theorem.

\begin{thm}
	For any nonnegative integers $n, k$, with $n \geq k$, we have
	\begin{equation}\label{Thm1}
	T_{2,\lambda}(n,k \mid x) = \sum_{l=k}^{n} \binom{n}{l} T_{2,\lambda}(l,k) (x)_{n-l,\lambda}.
	\end{equation}
\end{thm}

\vspace{0.1in}

Now, we observe that
\begin{equation}\label{19}
\begin{split}
& \frac{1}{k!} (1+\lambda t)^{\frac{x}{\lambda}} ((1+\lambda t)^{\frac{1}{2 \lambda}} - (1+\lambda t)^{-\frac{1}{2 \lambda}})^k \\
&= \frac{1}{k!} (1+\lambda t)^{\frac{1}{\lambda}(x-\frac{k}{2})} \sum_{l=0}^k \binom{k}{l} (-1)^{k-l} (1+\lambda t)^{\frac{l}{\lambda}} \\
& = \frac{1}{k!} \sum_{l=0}^k \binom{k}{l} (-1)^{k-l} e^{\frac{1}{\lambda}(x-\frac{k}{2}+l) \log(1+\lambda t)} \\
& = \frac{1}{k!} \sum_{l=0}^k \binom{k}{l} (-1)^{k-l} \sum_{m=0}^{\infty} \bigg(\frac{x-\frac{k}{2}+l}{\lambda} \bigg)^m \frac{1}{m!} \bigg(\log(1+\lambda t) \bigg)^m \\
& = \sum_{n=0}^{\infty} \frac{1}{k!} \sum_{l=0}^k \binom{k}{l} (-1)^{k-l} \sum_{m=0}^{n} \lambda^{n-m} S_1(n,m) (x+l-\frac{k}{2})^m \frac{t^n}{n!} \\
& = \sum_{n=0}^{\infty} \bigg( \sum_{m=0}^{n}(\frac{1}{k!} \delta^k x^{m}) \lambda^{n-m} S_1(n,m) \bigg).
\end{split}
\end{equation}

\noindent Therefore, by \eqref{16} and \eqref{19}, we obtain the following theorem.
\begin{thm}
	For any nonnegative integers $n,k $, we have
	\begin{equation}\label{Thm2}
	\begin{split}
	\sum_{m=0}^{n} (\frac{1}{k!} \delta^k x^{m}) \lambda^{n-m} S_1(n,m) =
	\begin{cases}
	T_{2,\lambda}(n,k \mid x),&\text{if}\,\, n \geq k,\\
	0,&\text{if}\,\, n < k.
	\end{cases}
	\end{split}
	\end{equation}
\end{thm}

\vspace{0.1in}
Letting $x=0$ in \eqref{Thm2} gives the next result.
\begin{thm}
	For any nonnegative integers $n,k $, we have
	\begin{equation}\label{Thm3}
	\begin{split}
	\sum_{m=0}^{n} (\frac{1}{k!} \delta^k 0^{m}) \lambda^{n-m} S_1(n,m) =
	\begin{cases}
	T_{2,\lambda}(n,k),&\text{if}\,\, n \geq k,\\
	0,&\text{if}\,\, n < k.
	\end{cases}
	\end{split}
	\end{equation}
\end{thm}

\vspace{0.1in}

By making use of \eqref{Thm2} and \eqref{02}, we note that
\begin{equation}\label{20}
\begin{split}
T_{2,\lambda}(n+1,k \mid x) & = \sum_{m=0}^{n+1} \bigg(\frac{1}{k!} \delta^k x^m \bigg) \lambda^{n+1-m} S_1(n+1,m) \\
& = \sum_{m=1}^{n+1} \bigg(\frac{1}{k!} \delta^k x^m \bigg) \lambda^{n+1-m} \bigg(S_1(n,m-1) - n S_1(n,m) \bigg) \\
& = \sum_{m=0}^{n} \bigg(\frac{1}{k!} \delta^k x^{m+1} \bigg) \lambda^{n-m} S_1(n,m) - n \lambda T_{2,\lambda}(n,k \mid x).
\end{split}
\end{equation}

\noindent On the other hand, by \eqref{15}, we get
\begin{equation}\label{21}
\begin{split}
& \sum_{m=0}^{n} \bigg(\frac{1}{k!} \delta^k x^{m+1} \bigg) \lambda^{n-m} S_1(n,m)\\
& = \sum_{m=0}^{n} \frac{1}{k!} \bigg( (x+\frac{k}{2}) \delta^k x^{m} + k \delta^{k-1} (x-\frac{1}{2})^{m} \bigg) \lambda^{n-m} S_1(n,m) \\
& = (x+\frac{k}{2}) \sum_{m=0}^{n}  \bigg( \frac{1}{k!} \delta^k x^m \bigg) \lambda^{n-m} S_1(n,m)\\
&\quad\quad +  \sum_{m=0}^{n} \bigg( \frac{1}{(k-1)!} \delta^{k-1}(x-\frac{1}{2})^m \bigg) \lambda^{n-m} S_1(n,m) \\
& = \bigg(x+\frac{k}{2}\bigg) T_{2,\lambda}(n,k \mid x) + T_{2,\lambda}(n,k-1 \mid x-\frac{1}{2}).
\end{split}
\end{equation}

\vspace{0.1in}

\noindent Therefore, by \eqref{20} and \eqref{21}, we obtain the following theorem.
\begin{thm}
	For any integers $n,k$, with $1 \leq k \leq n$, we have
	\begin{equation}\label{Thm4}
	T_{2,\lambda}(n+1,k \mid x) = \bigg(x+\frac{k}{2}-n \lambda\bigg) T_{2,\lambda}(n,k \mid x) + T_{2,\lambda}(n,k-1 \mid x-\frac{1}{2}).
	\end{equation}
\end{thm}

Setting $x=0$ in \eqref{Thm4} yields the following result.
\begin{thm}
	For any integers $n,k$, with $1 \leq k \leq n$, we have
	\begin{equation}\label{Thm5}
	T_{2,\lambda}(n+1,k) = \bigg(\frac{k}{2}-n \lambda\bigg) T_{2,\lambda}(n,k) + T_{2,\lambda}(n,k-1 \mid -\frac{1}{2}).
	\end{equation}
\end{thm}

Note that taking $\lambda \rightarrow 0$ in \eqref{Thm5} gives us
\begin{equation*}
	T(n+1,k) = \frac{k}{2} T(n,k) + T(n,k-1|-\frac{1}{2}),
\end{equation*}
where $1 \leq k \leq n$.

By \eqref{17}, we get
\begin{equation}\label{22}
\begin{split}
\frac{1}{k!} \bigg((1+\lambda t)^{\frac{1}{2 \lambda}} - (1+\lambda t)^{-\frac{1}{2 \lambda}}\bigg)^{k}
& = \frac{1}{k!} (1+\lambda t)^{-\frac{k}{2 \lambda}} \bigg((1+\lambda t)^{\frac{1}{\lambda}} -1 \bigg)^k \\
& = \frac{1}{k!} \sum_{l=0}^k \binom{k}{l} (-1)^{k-l} (1+\lambda t)^{\frac{1}{\lambda} (l-\frac{k}{2})} \\
& = \sum_{n=0}^\infty \bigg(\frac{n!}{k!} \sum_{l=0}^k \binom{k}{l} \binom{l-\frac{k}{2}}{n}_{\lambda} (-1)^{k-l} \bigg) \frac{t^n}{n!}.
\end{split}
\end{equation}

\noindent Therefore, by \eqref{17} and \eqref{22}, we obtain the following theorem.

\begin{thm}
	For any nonnegative integers $n,k $, we have
	\begin{equation*}\label{Thm6}
		\begin{split}
			\frac{n!}{k!} \sum_{l=0}^k \binom{k}{l} \binom{l-\frac{k}{2}}{n}_{\lambda} (-1)^{k-l} =
			\begin{cases}
				T_{2,\lambda}(n,k),&\text{if}\,\, n \geq k,\\
				0,&\text{if}\,\, n < k.
			\end{cases}
		\end{split}
	\end{equation*}
\end{thm}

\vspace{0.1in}

The Carlitz degenerate Euler polynomials of  higher-order (see\,\,\eqref{03}) are defined as
\begin{equation} \label{23}
\bigg(\frac{2}{(1+\lambda t)^{\frac{1}{\lambda}} +1} \bigg)^r (1+\lambda t)^{\frac{x}{\lambda}} = \sum_{n=0}^\infty \mathcal{E}_{n,\lambda}^{(r)} (x) \frac{t^n}{n!},~ ( r \in \mathbb{R}).
\end{equation}

\noindent Then \eqref{23} is also given by
\begin{equation}\label{24}
\begin{split}
& \bigg(\frac{2}{(1+\lambda t)^{\frac{1}{\lambda}} +1} \bigg)^r (1+\lambda t)^{\frac{x}{\lambda}} = 2^r \bigg((1+\lambda t)^{\frac{1}{\lambda}} +1 \bigg)^{-r} (1+\lambda t)^{\frac{x}{\lambda}} \\
& =  \bigg(\frac{(1+\lambda t)^{\frac{1}{\lambda}} -1} {2} +1 \bigg)^{-r} (1+\lambda t)^{\frac{x}{\lambda}} \\
& =  \sum_{l=0}^\infty \binom{r+l-1}{l} (-\frac{1}{2})^l ((1+\lambda t)^{\frac{1}{\lambda}} -1)^l (1+\lambda t)^{\frac{x}{\lambda}} \\
& =  \sum_{l=0}^\infty \binom{r+l-1}{l} (-\frac{1}{2})^l ((1+\lambda t)^{\frac{1}{2 \lambda}} - (1+\lambda t)^{-\frac{1}{2 \lambda}})^l (1+\lambda t)^{\frac{1}{\lambda}(x+\frac{l}{2})} \\
& =  \sum_{l=0}^\infty \binom{r+l-1}{l} (-\frac{1}{2})^l l! \sum_{n=l}^{\infty}  T_{2,\lambda}(n,l \mid x + \frac{l}{2}) \frac{t^n}{n!} \\
& =  \sum_{n=0}^\infty \bigg(\sum_{l=0}^n \binom{r+l-1}{l} (-\frac{1}{2})^l l! T_{2,\lambda}(n,l \mid x + \frac{l}{2}) \bigg) \frac{t^n}{n!}. \\
\end{split}
\end{equation}

\noindent Therefore, by \eqref{23} and \eqref{24}, we obtain the following theorem.
\begin{thm}
	For any nonnegative integer $n$, we have
	\begin{equation*}\label{Thm7}
		\mathcal{E}_{n,\lambda}^{(r)} (x) = \sum_{l=0}^n \binom{r+l-1}{l} (-\frac{1}{2})^l l!T_{2,\lambda}(n,l \mid x + \frac{l}{2}).
	\end{equation*}
\end{thm}

\vspace{0.1in}

By \eqref{17} and \eqref{05}, we get
\begin{equation}\label{25}
\begin{split}
& \sum_{n=2k}^\infty T_{2,\lambda}(n,2k) \frac{t^n}{n!} =  \frac{1}{(2k)!} \bigg((1+\lambda t)^{\frac{1}{2 \lambda}} - (1+\lambda t)^{-\frac{1}{2 \lambda}}\bigg)^{2k} \\
& =  \frac{1}{(2k)!} \bigg((1+\lambda t)^{\frac{1}{\lambda}} + (1+\lambda t)^{-\frac{1}{\lambda}} - 2 \bigg)^{k} \\
& =  \frac{1}{(2k)!} \sum_{l=0}^k \binom{k}{l} \bigg((1+\lambda t)^{\frac{1}{\lambda}} -1 \bigg)^l \bigg((1+\lambda t)^{-\frac{1}{\lambda}} -1 \bigg)^{k-l} \\
& =  \frac{k!}{(2k)!} \sum_{l=0}^k \sum_{n=k}^\infty \bigg(\sum_{i=l}^n \binom{n}{i} S_{2,\lambda}(i,l) S_{2,-\lambda}(n-i,k-l) (-1)^{n-i} \bigg) \frac{t^n}{n!} \\
& =  \sum_{n=k}^\infty \bigg(\frac{1}{k! \binom{2k}{k}} \sum_{l=0}^k \sum_{i=l}^n \binom{n}{i} S_{2,\lambda}(i,l) S_{2,-\lambda}(n-i,k-l) (-1)^{n-i} \bigg) \frac{t^n}{n!} .
\end{split}
\end{equation}

\noindent Comparing the coefficients on both sides of \eqref{25}, we have the following theorem.
\begin{thm}
	For any nonnegative integers $n,k $, we have
	\begin{equation*}\label{Thm8}
		\begin{split}
			\sum_{l=0}^k \sum_{i=l}^n \binom{n}{i} S_{2,\lambda}(i,l) S_{2,-\lambda}(n-i,k-l) (-1)^{n-i} =
			\begin{cases}
				k! \binom{2k}{k} T_{2,\lambda}(n,2k), &\text{if}\,\, n \geq 2k,\\
				0,&\text{if}\,\, n < 2k.
			\end{cases}
		\end{split}
	\end{equation*}
\end{thm}

Finally, we would like to define the degenerate central factorial numbers of the first kind.
We fisrt recall that the $\lambda$-logarithmic function is defined by
\begin{equation}\label{26}
\log_{\lambda} t = \frac{t^{\lambda} -1}{\lambda},
\end{equation}
where we note that $\lim_{\lambda \rightarrow 0} \log_{\lambda} t = \log t$.
Let $g(t)=(1+\lambda t)^{\frac{1}{2\lambda}}-(1+\lambda t)^{-\frac{1}{2\lambda}}$.

Then, using \eqref{26} we see that the inverse $g^{-1}(t)$ of $g(t)$ is given by
\begin{equation}\label{27}
g^{-1}(t)= \log_{\lambda} \bigg(\frac{t}{2}+ \sqrt{1+\frac{t^2}{4}}\, \bigg)^2. \\
\end{equation}

\noindent In view of \eqref{17} and \eqref{27}, we are led to define the degenerate central factorial numbers of the first kind by
\begin{equation*}
	\frac{1}{k!} \bigg(\log_{\lambda} \bigg(\frac{t}{2}+ \sqrt{1+\frac{t^2}{4}}\, \bigg)^{2}\bigg)^k = \sum_{n=k}^\infty t_{1,\lambda}(n,k) \frac{t^n}{n!}.
\end{equation*}

\bigskip
\medskip


\vspace{1cc}

{\bf{Taekyun Kim}} : Corresponding Author \\
Department of Mathematics\\
Kwangwoon University, Seoul, 139-701, Republic of Korea \\
Email:tkkim@kw.ac.kr \\

{\bf{Dae San Kim}}\\
Department of Mathematics \\
Sogang University, Seoul, 121-742, Republic of Korea \\
Email:dskim@sogang.ac.kr

{\small
	\noindent

}\end{document}